\documentclass[11pt,a4paper,twoside]{amsart}

\usepackage{amsfonts, amssymb,indentfirst}

\newtheorem{prop}{Proposition}[section]
\newtheorem{theorem}[prop]{Theorem}
\newtheorem{lemma}[prop]{Lemma}

\newtheorem{coro}[prop]{Corollary}
\newtheorem{remark}[prop]{Remark}
\newtheorem{con}[prop]{Conjecture}

\newcommand{\cqd}{\hfill$\Box$}
\newcommand{\mor}[0]{\operatorname{Mor}}

\newcommand{\ch}[0]{\operatorname{ch}}
\newcommand{\Td}[0]{\operatorname{Td}}
\newcommand{\rank}[0]{\operatorname{rank}}

\title[ON THE COHOMOLOGY OF BRILL-NOETHER LOCI OVER QUOT SCHEMES] {ON THE COHOMOLOGY OF BRILL-NOETHER LOCI OVER QUOT SCHEMES}
\author[Cristina Mart{\'\i}nez]{Cristina Mart{\'\i}nez}

\subjclass[2000]{ Primary 14F05; Secondary 14C40.} \keywords{Quot
schemes, Brill-Noether theory, vector bundles.} 
\address{Institut for Matematike Fag (CTQM), Aarhus Universitet,
Ny Munkegade, 8000 Aarhus, Denmark} \email{cmartine@imf.au.dk}
\address{Max Planck Institute for Mathematics, Vivatgasse 7, Bonn, 53111.}
 \email{cmartine@mpim-bonn.mpg.de}
\date{Jannuary 2008}
\begin{document}
\maketitle
\begin{abstract}
Let $C$ be a smooth projective irreducible curve over an algebraic
closed field $k$ of characteristic 0. We consider Brill-Noether
loci over the moduli space of morphisms from $C$ to a Grassmannian
$G(m,n)$ of $m-$planes in $k^{n}$ and the corresponding Quot
schemes of quotients of a trivial vector bundle on $C$
compactifying the spaces of morphisms. We study in detail the case
in which $m=2, n=4$. We prove results on the irreducibility and
dimension of these Brill-Noether loci and we address explicit
formulas for their cohomology classes. We study the existence
problem of these spaces which is closely related with the problem
of classification of vector bundles over curves.
\end{abstract}


\normalsize \setcounter{equation}{0}
\section{Introduction}
 The theory of Brill-Noether
over the space of stable vector bundles or semistable bundles has
been very much studied (see for example \cite{BGMN}, \cite{BGN},
\cite{Mer}, \cite{Te}). Let $C$ be a non-singular projective curve
defined over an algebraically closed field $k$ of characteristic
0, and let $M(r,d)$ denote the moduli space of stable vector
bundles over $C$ of rank $r$ and degree $d$. When $r$ and $d$ are
not coprime, $M(r,d)$ is not compact and can be compactified to a
scheme $\widetilde{M}(r,d)$ by adding equivalence classes of
semi-stable vector bundles of rank $r$ and degree $d$. For $E$ a
rank $r$ and degree $d$ vector bundle, the slope of $E$ is defined
as $\mu\,(E)=\frac{d}{r}$. The notion of stability, semistability
and $S-$equivalence was first introduced by Mumford, Seshadri and
Narasimhan, \cite{Mum}, \cite{NS}, \cite{Ses}. A vector bundle is
stable (respectively semistable) when for every proper subbundle
$F$;
$$\mu(F)=\frac{deg\,(F)}{rank\,(F)}< \mu(E), \ \rm{respectively}  \leq.$$

Any semistable vector bundle $E$ has a $S-$filtration, that is, a
filtration by subbundles:
$$0=E_{0}\subset E_{1} \subset \ldots \subset E_{k}=E, $$

\noindent whose factors $(E_{i}/E_{i-1})$ are all stable with $\mu
(E_{i}/E_{i-1})=\mu(E)$. The isomorphism class of the direct sum
$gr\,E:=\bigoplus_{i=1}^{k} E_{i}/E_{i-1}$ is independent of the
filtration, and two semistable bundles are called $S-$equivalent
if
$$gr\,(E)\cong gr\,(F).$$


The Brill-Noether loci over the moduli space of stable bundles are
defined as:

$$B_{r,d,k}=\{E\in M(r,d)|\, h^{0}(E)\geq k \}$$

\noindent for a fixed integer $k$, and over the moduli space of
semistable vector bundles it is defined by:

$$\widetilde{B}_{r,d,k}=\{[E]\in M(r,d)| h^{0}(gr\,(E))\geq k\}.$$



\noindent By the Semicontinuity theorem, these Brill-Noether loci
are closed subschemes of the appropriate moduli spaces, and in
particular it is not difficult to describe them as determinantal
loci which allows us to estimate their dimension.

The main object of Brill-Noether theory is the study of these
subschemes, in particular  questions related to their
non-emptiness, connectedness, irreducibility, dimension,
topological and geometric structure. In the case of line bundles
in which the moduli spaces are all isomorphic to the Jacobian,
these questions have been completely answered when the underlying
curve is generic (see for instance \cite{ACGH}). Montserrat
Teixidor i Bigas  proved in 1991 a theorem that gives a rather
general solution to the problem, but it holds only for a generic
curve (see \cite{Te}). Brambila-Paz, Grezegorczyk and Newstead in
1995 gave a solution to the problem for every curve when
$\frac{d}{r}\leq 1$, (cf \cite{BGN}).

We can define in an analogous way Brill-Noether loci over the
space of morphisms  $\mor_{d}\,(C,X)$ of fixed degree $d$ from a
curve to a projective variety, and the corresponding Quot schemes
compactifying these spaces of morphisms. In particular, we are
going to consider a Brill-Noether stratification over the space of
morphisms $R^{0}_{C,d}$ from a genus $g$, smooth projective
irreducible curve $C$ to the Grassmannian $G(m,n)$ of $m$
dimensional subspaces of $k^{n}$. Moreover, we study in detail the
case in which $m=2$, $n=4$. In this case, this theory is connected
with the geometry of ruled surfaces in $\mathbb{P}^{3}$. Since
$G(2,4)$ is endowed with the universal quotient bundle
$\mathcal{Q}$, a natural compactification for this space is
provided by the Quot scheme $R_{C,d}$, parametrizing epimorphisms
of vector bundles $\mathcal{O}^{4}_{C}\rightarrow E\rightarrow 0$.
Quot schemes have been shown by Grothendieck, to be fine moduli
spaces for the problem of parametrizing quotients of a fixed
sheaf, and as such, to carry universal structures.

Given a morphism $f$ in $R^{0}_{C,d}$, the pull-back
$E:=f^{*}\mathcal{Q}^{\vee}$ is a vector bundle of rank 2 over
$C$. We can consider its Segre invariant $s(E)$ which is defined
as the minimal degree of a twist $E^{\vee}\otimes L$ with a line
bundle such that the resulting bundle has a non-zero section. The
Brill-Noether stratum is then defined by the following locally
closed condition:
$$R^{0}_{C,d,s}:=\{f\in R^{0}_{C,d}|\,{\bf
s}(f^{*}\mathcal{Q})=s\}.
$$

We consider the closure of the stratum $R_{C,d,s}$ in the Quot
scheme compactification. In the first three sections of the paper
some basic results on the loci $R_{C,d,s}$ are presented. In
particular $R_{C,d,s}$ is exhibited as the degeneration locus of a
natural and appropriate morphism of vector bundles.

In section four the fundamental class of $R_{C,d,s}$ is computed
in the cohomology ring of $R_{C,d}$ under the assumption that
$R_{C,d,s}$ is either empty or of the expected codimension in
$R_{C,d}$ (for large $d$ depending on $s$), (Theorem \ref{class}).
We give a
partial solution to the existence problem in the rank two case,
that is, when we are considering the Grassmannian $G(2,4)$. We see
that the problem is closely related with the problem of
classification of vector bundles over curves.






\section{Brill-Noether loci}

Consider the universal exact sequence over the Grassmannian
$G(m,n)$:

\begin{equation}\label{seq1} 0\rightarrow \mathcal{N}\rightarrow
\mathcal{O}^{n}_{G}\rightarrow \mathcal{Q}\rightarrow 0\
\end{equation}

\noindent For every morphism $f\in M_{d}:=\mor_{d}\,(C,G(m,n))$,
we take the pull-back of the sequence \ref{seq1}:

\begin{equation}
0\rightarrow f^{*}\mathcal{N}\rightarrow
f^{*}\mathcal{O}^{n}_{G(m,n)}\rightarrow
f^{*}\mathcal{Q}\rightarrow 0.
\end{equation}




The next lemma describes the bundles $E$ which arise from this
construction.



\begin{lemma}\label{sections} Given $E$, a degree $d$, rank $r$ bundle over $C$,
there exists a morphism $f\in \mor_{d}(C,G(n-r,n))$ such that
$f^{*}\mathcal{Q}=E$ if and only if $E$ is generated by $n$ global
sections or equivalently is given by a quotient,
$$\mathcal{O}^{n}_{C}\rightarrow E\rightarrow 0.$$
\end{lemma}


\noindent We consider the Grassmannian $G(m,n)$ where $m=n-r$. By
the universal property of the Grassmannian, there exists a
morphism $f\in M_{d}$ such that $f^{*}\mathcal{Q}\cong E$, where
$\mathcal{Q}$ is the universal quotient bundle over the
Grassmannian $G(m,n)$.

\noindent Conversely, for all $f\in M_{d}$, $f^{*}\mathcal{Q}$ is
generated by global sections, since $\mathcal{Q}$ is given by a
quotient:

$$\mathcal{O}^{n}_{G(m,n)}\rightarrow \mathcal{Q}\rightarrow 0.$$

These quotients are parametrized by Grothendieck's Quot schemes
$Q_{d,r,n}(C)$ of degree $d$, rank $r$ quotients of
$\mathcal{O}_{C}^{n}$ compactifying the spaces of morphisms
$M_{d}$.

In the genus 0 case, by a theorem of Grothendieck, every vector
bundle $E$ over $\mathbb{P}^{1}$ decomposes as a direct sum of
line bundles and therefore for $E\cong \bigoplus_{i}
\mathcal{O}_{\mathbb{P}^{1}}(a_{i})$ to be generated by global
sections means that $a_{i}\geq 0$. In the genus 1 case, the
bundles generated by global sections are the indecomposable
bundles of degree $d>r$ (since $h^{0}(E)=d>r$ for an
indecomposable bundle of positive degree \cite{At}), the trivial
bundle $\mathcal{O}_{C}$ and the direct sums of such bundles. In
genus greater than or equal to 2, certain restrictions on the
bundle imply that it is generated by global sections. For example,
we can tensorize with a line bundle of degree $m$ such that
$E(m)=E\otimes \mathcal{O}_{C}(m)$ is generated by global sections
and $h^{1}(E(m))=0$, \cite{Ses}. Moreover, if $d$ is sufficiently
large and $E$ is semistable, then $h^{1}(E)=0$ and $E$ is
generated by its global sections; in fact, it is sufficient to
take $d>r\,(2g-1)$.
\vspace{0.5cm}

\noindent We {\it define the Brill-Noether loci over the spaces of
morphisms $M_{d}$} as:
\begin{equation}\label{BNL}
M_{d,a}=\{f\in M_{d}|\,h^{0}(C,f^{*}\mathcal{Q})\geq a\}
\end{equation}
for a fixed integer $a$. More generally, we can tensorize the
bundle with a fixed line bundle $L$ over $C$ and consider the
following Brill-Noether loci:

\begin{equation}\label{BNLL}
M_{d,a}(L)=\{f\in M_{d}|\, h^{0}(C,f^{*}\mathcal{Q}\otimes L)\geq
a\}.
\end{equation}

\section{A Brill-Noether stratification over the Quot scheme.}
In \cite{Mar2}, we consider the space of morphisms
$R^{0}_{d}:=\mor_{d}(\mathbb{P}^{1},G(2,4))$ and two different
compactifications of this space, the Quot scheme compactification
and the compactification of stable maps given by Kontsevich. We
consider the following Brill-Noether loci inside the space of
morphisms $R^{0}_{d}$:

\begin{equation}\label{Rd,a}
R^{0}_{d,a}=\{f\in R^{0}_{d}|\,
h^{0}\,(f^{*}\mathcal{Q}^{\vee}\otimes
\mathcal{O}_{\mathbb{P}^{1}}(a))\geq 1, \
h^{0}\,(f^{*}\mathcal{Q}^{\vee}\otimes
\mathcal{O}_{\mathbb{P}^{1}}(a-1))=0 \},
\end{equation}
\vspace{0.5cm} \noindent for a fixed integer $a$.


\noindent Note that we are considering here rank two bundles, but
the definition can be generalized easily to bundles of arbitrary
rank $r$.

It is easy to see that this set can be defined alternatively as
the $f\in R^{0}_{d}$ with $f^{*}\mathcal{Q}\cong
\mathcal{O}_{\mathbb{P}^{1}}(a)\oplus
\mathcal{O}_{\mathbb{P}^{1}}(d-a)$, for $a\leq \frac{d}{2}$, and
the parameter $a$ gives a stratification of the space $R^{0}_{d}$.

\subsubsection*{Geometric interpretation} The image of a curve $C$
by $f$ is a geometric curve in the corresponding Grassmannian or
equivalently a rational ruled surface in $\mathbb{P}^{3}$ for the
Grassmannian of lines. Fixing the parameter $a$ we are fixing the
degree of a minimal directrix in the ruled surface. The spaces
$R^{0}_{d,a}$ are locally closed again by the Semicontinuity
Theorem and they can be shown as the degeneration locus of a
morphism of bundles by means of the universal exact sequence over
the corresponding Quot scheme $R_{d}$ and we find that the
expected dimension of $R^{0}_{d,a}$ as a determinantal variety is
$3d+2a+5$. These spaces are considered in \cite{Mar2} as parameter
spaces for rational ruled surfaces in order to solve the following
enumerative problems:
\begin{enumerate}
\item The problem of enumerating rational ruled surfaces through
$4d+1$ points, or equivalently computing the degree of $R^{0}_{d}$
inside the projective space of surfaces of fixed degree $d$,
$$R^{0}_{d}\rightarrow \mathbb{P}^{{d+3 \choose 3}-1}.$$ 

\item Enumerating rational ruled surfaces with fixed splitting
type. This problem raises the question of defining Gromov-Witten
invariants for bundles with a fixed splitting type. 
\end{enumerate}

\vspace{0.5cm} From now, the underlying curve $C$ will be a
smooth, irreducible projective curve of genus greater than or
equal to 1, and we will be studying only the case $n=2, m=4$. We
call $R^{0}_{C,d}$ the spaces of morphisms $\mor_{d}\,(C,G(2,4))$.
For a vector bundle $E$ of rank 2 over $C$, its Segre invariant is
the integer $s$ such that the minimal degree of a line quotient
$E\rightarrow L \rightarrow 0$ is $\frac{d+s}{2}$, or the maximal
degree of a twist $E^{\vee}\otimes L$, having a nonzero section.
 Note that
$s(E)\equiv deg(E) \ {\rm{mod}} \ 2$ and that $E$ is stable (resp.
semistable)
if and only if $s(E)\geq 1$ (resp. $s(E)\geq 0$), (\cite{LN}, \cite{CS}).

If $T$ is any algebraic variety over $k$ and $\mathcal{A}$ is a
vector bundle of rank $r$ on $C\times T$, then the function
$s:T\rightarrow \mathbb{Z}$ defined by
$s(t)=s(\mathcal{A}|_{C\times t})$ is lower semicontinuous. 

Given $f\in R^{0}_{C,d}$ we consider the Segre invariant $s$ of
the bundle $f^{*}\mathcal{Q}$ which is the maximal degree of a
twist $f^{*}\mathcal{Q}^{\vee}\otimes L$ with a generic line
bundle $L$ of degree $\frac{d+s}{2}$ such that the resulting
bundle has a non-zero section.

We {\it define the corresponding Brill-Noether loci over
$R^{0}_{C,d}$} as the subsets:

\begin{equation}\label{BNL}
R^{0}_{C,d,s}
 =\{f\in
R^{0}_{C,d}|\,h^{0}(f^{*}\mathcal{Q^{\vee}}\otimes L)\geq 1, \ \
 L {\rm{\ of \ minimal \ degree}} \ \frac{d+s}{2}\}.
 \end{equation}

 The subvarieties $R^{0}_{C,d,s}$ are locally closed and are a
natural generalization for arbitrary genus of the subvarieties
$R^{0}_{d,a}$, defined previously, taking $a=\frac{d+s}{2}$ which
is the degree of the line bundle $L$.

The Zariski closure of $R^{0}_{C,d,s}$ is given by the set,
(\cite{Mar2}):
$$\overline{R}^{0}_{C,d,s}=\{f\in R^{0}_{C,d}| h^{0}(f^{*}\mathcal{Q}^{\vee}
\otimes L)\geq 1, \ \ deg(L)=\frac{d+s}{2}\}.$$

Observe that $L$ is allowed to vary on the variety of line bundles
of fixed degree, which is isomorphic to the Jacobian, which is
compact. The closure describes exactly the set of morphisms $f\in
R^{0}_{C,d}$ for which $f^{*}\mathcal{Q}$ is a rank two vector
bundle with Segre invariant less or equal than $s$, that is,

$$\overline{R}^{0}_{C,d,s}=R^{0}_{C,d,s}\cup R^{0}_{C,d,s-1} \cup \ldots $$

Evenmore, $R^{0}_{C,d,s-1}\subset R^{0}_{C,d,s}$, and therefore
$\overline{R}^{0}_{C,d,s-1}\subset \overline{R}^{0}_{C,d,s}$.

Let us consider the universal exact sequence on $R_{C,d}\times C$,

\begin{equation}\label{us}
0\rightarrow \mathcal{K}\rightarrow \mathcal{O}^{n}_{R_{C,d}\times
C}\rightarrow \mathcal{E}\rightarrow 0,
\end{equation}
which satisfies that for all $k-$scheme $S$, the set of morphisms
$f:S\rightarrow R_{C,d}$ is in correspondence one to one with the
set of isomorphism classes of short exact sequences,
$\mathcal{O}^{4}_{S\times C}\rightarrow \mathcal{E}_{S\times
C}\rightarrow 0$, where $\mathcal{E}_{S\times C}$ is flat over
$S$.

 Since the universal quotient sheaf $\mathcal{E}$ is flat over
$R_{C,d}$, for each $q\in R_{C,d}$,
$E_{q}:=\mathcal{E}\mid_{\{q\}\times C}$ is a coherent sheaf over
$C$. By flatness, $h^{0}(E_{q})-h^{1}(E_{q})$ is constant on any
connected component of  $R_{C,d}$. The Riemman-Roch theorem allows
us to compute its value:
$$h^{0}(E_{q})-h^{1}(E_{q})=d+2\,(1-g), $$ for every
$q\in R_{C,d}$.

We consider the Zariski closure $R_{C,d,s}$ of the sets
$R^{0}_{C,d,s}$ inside the Quot scheme compactification of the
space of morphisms:

$$R_{C,d,s}=\{q\in R_{C,d}|\,h^{0}\,(C,E_{q}^{\vee}\otimes L)\geq 1,
\ deg\,L=\frac{d+s}{2}\}$$ $$=\{q\in
R_{C,d}|\,h^{1}\,(C,E_{q}\otimes  K_{C}\otimes L^{-1})\geq 1, \
deg\,L=\frac{d+s}{2}\}$$
$$=\{q\in R_{C,d}|\,h^{0}\,(C,E_{q}\otimes K_{C} \otimes  L^{-1})\geq d+3-2g,
\ deg\,L=\frac{d+s}{2}\}.$$ The first equality is due to the
duality theorem, and the second one is due to
 the Riemann-Roch theorem.

\subsubsection*{The existence problem}
The existence problem here means that, given $s$, there exist a
vector bundle $E$ of rank 2 with Segre invariant $s$ and a
morphism $f\in R^{0}_{C,d}$ such that $f^{*}\mathcal{Q}:=E$.
Therefore the existence problem for $R^{0}_{C,d,s}$ is closely
related with the existence problem of vector bundles of rank 2
with Segre invariant $s$ and consequently with the problem of
classification of vector bundles over curves. The contents of the
next proposition are well known in the context of vector bundles,
(see for example \cite{La} for the rank 2 case and \cite{RT} for
the general case).

\begin{prop}\begin{enumerate}\item If $s>g$ then $R^{0}_{C,d,s}$ is empty.
\item If $C$ is a smooth elliptic curve, $d\equiv\, s\ mod\ 2$ and
$d\geq 3$ then $R^{0}_{C,d,1}$ and $R^{0}_{C,d,0}$ are non-empty.
\item If $C$ is a smooth curve of genus $g\geq s \geq 0$,
$d\equiv\, s\ mod\ 2$ and $d>2\,(2g-1)$, then $R^{0}_{C,d,s}$ is
non empty.
\end{enumerate}
\end{prop}
{\it Proof.} Given $f\in R^{0}_{C,d}$, the pull-back
$f^{*}\mathcal{Q}^{\vee}$ of the dual of the universal quotient
bundle over $G(2,n)$ is a rank 2 bundle over $C$ of degree $-d$.
We can tensorize it with a line bundle of degree $\frac{d+s}{2}$
such that $E':=f^{*}\mathcal{Q}^{\vee}\otimes L$ is a normalized
bundle of degree $s$, that is, $h^{0}(E')\neq 0$ but for all
invertible sheaves $\mathcal{L}$ on $C$ with $deg\,\mathcal{L}<0$,
we have $H^{0}(E'\otimes L)=0$. Since $E'$ is a normalized bundle
there are two possibilities for $E'$:
\begin{itemize}
\item If $E'$ is decomposable as a direct sum of two invertible
sheaves, then $E'\cong \mathcal{O}_{C}\oplus F$ for some $F$ with
$deg F\leq 0$, therefore $s(E')\leq 0$ and all values of $s(E)\leq
0$ are possible.

\item If $E'$ is indecomposable, then $2-2g\leq s\leq g$, (see
\cite{Har}).
\end{itemize}
In particular, this implies that there does not exist a vector
bundle of rank 2 over $C$ with Segre invariant $s>g$, and
consequently $R^{0}_{C,d,s}$ is empty as we stated in
(1).

If $C$ is an elliptic curve, for each value $0\leq s \leq 1$ there
is a bundle $E$ over $C$ of rank 2 with Segre invariant $s$.
Therefore by Lemma \ref{sections}, the sets $R^{0}_{C,d,1}$ and
$R^{0}_{C,d,0}$ are non-empty. This proves (2).

If $C$ is a smooth curve of genus $g\geq 2$, $0\leq s\leq g$ and
$d>2\,(2g-1)$, there is a semistable vector bundle of rank 2 over
$C$ with Segre invariant $s$ and again by Lemma \ref{sections},
$R^{0}_{C,d,s}$ is non-empty. \cqd

\vspace{0.5cm} The next theorem will exhibit the $R_{C,d,s}$ as
determinantal varieties which allows us to estimate their
dimensions.


Let $K_{C}$ be the canonical bundle over $C$ and $\pi_{1},
\pi_{2}$ be the projection maps of $R_{C,d}\times C$ over the
first and second factors respectively. Tensorizing the sequence
(\ref{us}) with the linear sheaf $\pi_{2}^{*}(K_{C}\otimes
L^{-1})$ gives the exact sequence:

\begin{equation}\label{seq}
0\rightarrow \mathcal{K}\otimes \pi_{2}^{*}(K_{C}\otimes
L^{-1})\rightarrow \mathcal{O}^{n}_{R_{C,d}\times C}\otimes
\pi_{2}^{*}(K_{C}\otimes L^{-1})\rightarrow \mathcal{E}\otimes
\pi_{2}^{*}(K_{C}\otimes L^{-1})\rightarrow 0.
\end{equation}

Here $L$ is a generic line bundle on $C$ of fixed degree
$\frac{d+s}{2}$. The $\pi_{1*}$ direct image of the above sequence
yields the following long exact sequence on $R_{C,d}$:

$$
0\rightarrow \pi_{1*}(\mathcal{K}\otimes \pi_{2}^{*}(K_{C}\otimes
L^{-1}))\rightarrow \pi_{1*}( \mathcal{O}^{n}_{R_{C,d}}\otimes
\pi_{2}^{*}(K_{C}\otimes L^{-1}))\rightarrow$$ $$\rightarrow
\pi_{1*}(\mathcal{E}\otimes \pi_{2}^{*}(K_{C}\otimes
L^{-1}))\rightarrow R^{1}\pi_{1*}(\mathcal{K}\otimes
\pi_{2}^{*}(K_{C}\otimes L^{-1}))\rightarrow $$ $$ \rightarrow
R^{1}\pi_{1*}( \mathcal{O}^{n}_{R_{C,d}}\otimes
\pi_{2}^{*}(K_{C}\otimes L^{-1})) \rightarrow
R^{1}\pi_{1*}(\mathcal{E}\otimes \pi_{2}^{*}(K_{C}\otimes
L^{-1}))\rightarrow 0.
$$

\begin{theorem}\label{teo1} For $d$ sufficiently large depending on $s$,
$R_{C,d,s}(L)$ is the locus where the map
\begin{equation}\label{conucleo}R^{1}\pi_{1*}(\mathcal{K}\otimes
\pi^{*}_{2}(L^{-1}\otimes K_{C})\rightarrow
R^{1}\pi_{1*}(\mathcal{O}^{n}_{R_{C,d}\times C }\otimes
\pi^{*}_{2}(L^{-1}\otimes K_{C}))\end{equation} is not surjective.
It is irreducible and has expected codimension $2g-s-1$ as a
determinantal variety.

\end{theorem}\label{irred}
{\it Proof.} The map (\ref{conucleo}) is not surjective in the
support of the sheaf
$$R^{1}\pi_{1*}(\mathcal{E}\otimes
\pi_{2}^{*}(K_{C}\otimes L^{-1})),$$ that is, in the points $q\in
R_{C,d}$ such that $h^{1}\,(\mathcal{E}\otimes
\pi_{2}^{*}(K_{C}\otimes L^{-1})|_{\{q\}\times C})\geq 1$, or
equivalently in $R_{C,d,s}(L)$ by Serre duality. In other words,
by semicontinuity there is an open set $P_{s}\subset
\mathbb{P}(R^{1}\pi_{1*}(\mathcal{E}\otimes
\pi_{2}^{*}(K_{C}\otimes L^{-1})))$ parametrizing  classes of
quotients $\mathcal{O}^{4}_{C}\rightarrow E_{s}\rightarrow 0$ such
that $E_{s}$ is a rank 2 bundle with Segre invariant $s$, modulo
the canonical operation of $k^{*}$.

Let $P_{s}\stackrel{f}{\rightarrow}R_{C,d,s}$ be the surjective
morphism such that the image by $f$ of each quotient is its class
of isomorphism in $R_{C,d,s}(L)$. The fibers are isomorphic to
$\mathbb{P}(Aut(E_{s}))$.
This proves that the subschemes $R_{C,d,s}(L)$, being the image of
an irreducible variety by a morphism, are irreducible.

\noindent By Serre Duality it follows that
$$\left(R^{1}\pi_{1*}\,(\mathcal{K}\otimes \pi_{2}^{*}(L^{-1}\otimes K_{C}))\right)^{\vee} \cong
\pi_{1*}(\mathcal{K}^{\vee}\otimes \pi_{2}^{*}L)$$ and by the Base
Change Theorem, their fibers are isomorphic to
$$H^{0}(C,\mathcal{K}^{\vee}\otimes \pi_{2}^{*}L|_{C\times \{p\}}),\ \ p\in R_{C,d}$$
and have dimension $2d+s+m\,(1-g)$, where
$m=\,rank(\mathcal{K}^{\vee})$. It is enough to take
$d+s>2m\,(g-1)$ to ensure the vanishing of
$h^{1}(C,\mathcal{K}^{\vee}\otimes \pi_{2}^{*}L|_{C\times
\{p\}})$. As a consequence, $\pi_{1*}(\mathcal{K}^{\vee}\otimes
\pi_{2}^{*}L)$ is a bundle of rank $2d+s+2\,(1-g)$, (note that we
are specializing $m$ and $n$ to be 2, 4 respectively). Again by
Serre Duality we see that
$R^{1}\pi_{1*}(\mathcal{O}^{n}_{R_{C,d}\times C }\otimes
\pi^{*}_{2}(L^{-1}\otimes K_{C}))\cong
\pi_{1*}(\mathcal{O}^{n}_{R_{C,d}}\otimes \pi_{2}^{*}L)$ and it is
a bundle with fiber isomorphic to
$$H^{0}(C, \mathcal{O}^{n}_{R_{C,d}}\otimes \pi_{2}^{*}L|_{C\times \{p\}})$$
of dimension $2d+2s-4g+4$. Therefore we have the following
morphism of bundles:
$$\pi_{1*}(\mathcal{K}^{\vee}\otimes \pi_{2}^{*}L)\stackrel{\phi}{\rightarrow}
\mathcal{O}^{2d+2s-4g+4}_{R_{C,d}}.$$ The expected codimension of
$R_{C,d,s}$ as determinantal variety is
$$\left((2d-2g+s+2)-(2d+2s-4g+3)\right)$$ $$ \cdot\left( (2d+2s-4g+4)-(2d+2s-4g+3))\right.=2g-s-1.$$
\cqd

\begin{remark}Note that as $s$ increases, the bundle $E$ becomes
more general, so the codimension decreases. In particular when
$s=g-1$, the codimension is $g$ for fixed $L$ or 0 when we allow
$L$ to vary.
\end{remark}

\section{Cohomology of the varieties $R_{C,d,s}$.}
In this section we determine the cohomology classes of the Brill-Noether loci in
terms of some natural elements in the cohomology ring of $R_{C,d}$ defined by
the universal bundle.

Let $\{1,\delta_{k}, 1\leq k\leq 2g, \eta\}$ be a basis for the
cohomology of $C$, where $\eta$ represents the class of a point.
We will also denote by $\{1,\delta_{k}, 1\leq k\leq 2g, \eta\}$
the
pull-backs to $R_{C,d}$ by the projection morphism.

\noindent Let us consider the classes $t_{i},u_{i-1}, s_{i}^{j}$
in $H^{*}(R_{C,d};\mathbb{Q})$ defined by the K\"unneth
decompostion of the Chern clases of $\mathcal{K}^{\vee}$:

\noindent
$$c_{i}(\mathcal{K}^{\vee})=t_{i}+\sum_{j=1}^{2g}s_{i}^{j}\delta_{j}+u_{i-1}\,
\eta,\ t_{i}\in H^{2i}\,(R_{C,d};\mathbb{Q}), s^{j}_{i}\in
H^{2i-1}(R_{C,d};\mathbb{Q}),$$ $$ \ u_{i-1}\in
H^{2i-2}\,(R_{C,d};\mathbb{Q}).$$

\noindent Every class $z$ in the cohomology ring
$H^{*}(R_{C,d};\mathbb{Q})$ can be written in the form
$$z=c+\sum_{j=1}^{2g}b^{j}\delta_{j}+f\,\eta$$
where $c=\pi_{*}(\eta\,z)$ and $f=\pi_{*}(z)\in
H^{*}(R_{C,d};\mathbb{Q})$. In particular,
$t_{i}=\pi_{*}(\eta\,c_{i}(\mathcal{K}^{\vee}))$ and
$u_{i-1}=\pi_{*}(c_{i}(\mathcal{K}^{\vee}))$,
$u_{0}=\pi_{*}(c_{1}(\mathcal{K}^{\vee}))=d$.

\begin{con}The elements $t_{1}, t_{2}, u_{1}, s^{j}_{i} \,(1\leq j\leq 2g, i=1,2)$
generate $H^{*}(R_{C,d};\mathbb{Q})$.
\end{con}

The evidence for the conjecture is that in the genus 0 case it is
true (\cite{Str}) and $s^{j}_{i} (1\leq j\leq 2g, i=1,2)$ are the
generators that appear when we consider curves in genus higher
than 0. In addition, for the ordinary moduli space $M(2,d)$ with
$d$ odd, the $s^{j}_{i}$  are all needed, together with $t_{2}$
and one of $t_{1}, u_{1}$ as well (see \cite{Za}). For $d$ even,
the conjecture is less plausible, and it is unlikely to be true
for integer cohomology.

\noindent For $L$ a line bundle over $C$ of degree
$a=\frac{d+s}{2}$, its first Chern class is given by

$$c_{1}(\pi_{2}^{*}L)=a\,\eta.$$

\begin{theorem}\label{class}If $C$ is any smooth curve of genus $g$, and
$R_{C,d,s}$ is either empty or generically reduced and of the
expected codimension $2g-s-1$, $R_{C,d,s}$ has fundamental class:

$$[R_{C,d,s}]=-c_{2g-s-1}(\pi_{1*}(\mathcal{K}^{\vee}\otimes \pi^{*}_{2}L)).$$
\end{theorem}
{\it Proof.} By assumption, $R_{C,d,s}$ has the expected
codimension as a determinantal variety and is irreducible as we
have seen in Theorem \ref{teo1}; therefore, it does not have
components contained at infinity and the Porteous formula gives
the fundamental class of the varieties $R_{C,d,s}$ in terms of the
Chern classes of the bundles given by Theorem \ref{teo1}. We get
that
$$[R_{C,d,s}]=\triangle_{2g-s-1,1}\left(c_{t}(-\pi_{1*}(\mathcal{K}^{\vee}\otimes
\pi_{2}^{*}L))\right.,$$ where

$\triangle_{p,q}(a)=det\,\left(\begin{array}{ccc}
a_{p}&\ldots&a_{p+q-1}\\ \vdots& &\vdots\\
a_{p-q+1}&\ldots & a_{p}
\end{array}\right),$ for any formal series
$a(t)=\sum_{k=-\infty}^{k=+\infty}a_{k}\,t^{k}$.

\begin{remark}Note that the assumption that $R_{C,d}$ is generically reduced and of expected
codimension $2g-s-1$ is satisfied for large $d$ (relative to $s$
and $g$). In that case, $R_{C,d}$ is a smooth projective bundle
over the Jacobian $J^{d}$ of degree $d$ line bundles, and the
intersection numbers on $Quot_{d,r,n}(C)$ correspond to certain
counts of maps from $C$ to $G(r,n)$. The intersection
 of the $t$ classes has been studied extensively in \cite{Ber}, where it is shown that the evaluation of a
 top-degree monomial in the $t$ classes on the fundamental cycle has enumerative meaning. It is the number
 of degree $d$ maps from $C$ to the Grassmannian $G(r,n)$ which sends fixed points on $C$ to special Schubert
 varieties of $G(r,n)$. 
\end{remark}

\subsubsection{Computations of Chern classes}
Applying Grothendieck-Riemann-Roch theorem to the projection
morphism $\pi_{1}$, it follows that
\begin{equation}\label{RRG}\ch\,(\pi_{1*}(\mathcal{K}^{\vee}\otimes
\pi_{2}^{*}L)))=\pi_{1*}(\Td\,(R_{C,d}\times C)/R_{C,d})\cdot
\ch\,(\mathcal{K}^{\vee}\otimes \pi_{2}^{*}L)).\end{equation}

First we compute the Chern classes of $\mathcal{K}^{\vee}\otimes
\pi_{2}^{*}L$:
$$c_{1}(\mathcal{K}^{\vee}\otimes
\pi_{2}^{*}L))=t_{1}+(\sum_{j=1}^{2g}s_{1}^{j}
\delta_{j})+(d+2a)\,\eta,$$

$$c_{2}\,(\mathcal{K}^{\vee}\otimes \pi_{2}^{*}L)=t_{2}+(\sum_{j=1}^{2g}s_{2}^{j}
\delta_{j})+u_{1}\eta +a\eta t_{1}+
 a(\sum_{j=1}^{2g}s_{1}^{j}\delta_{j})\eta, $$


 Let $\alpha_{1}$ and $\alpha_{2}$ be the
 classes $\sum_{j=1}^{2g}s_{1}^{j} \delta_{j}$ and $\sum_{j=1}^{2g}s_{2}^{j}
\delta_{j}$ respectively. The intersection numbers for the
$\delta_{i}$ imply the following relations:

$$\alpha_{1}^{2}=-2A\eta, \ \ \ A=\sum_{j=1}^{g}s_{1}^{j}s_{1}^{j+g}\in H^{2}\,(R_{C,d};\mathbb{Q}), \ \ \alpha_{1}^{3}=0,$$

$$ \alpha_{2}^{2}=-2\gamma \eta, \ \ \ \gamma=\sum_{i=1}^{g}s_{2}^{j}s_{2}^{j+g}\ \in H^{6}\,(R_{C,d};\mathbb{Q}),
\ \alpha_{2}^{3}=0,$$

$$\alpha_{1}\alpha_{2}=B\eta, \ \ B=(\sum_{i=1}^{g}-s_{1}^{i}s_{2}^{i+g}+s_{1}^{i+g}s_{2}^{i})\in H^{4}(R_{C,d};\mathbb{Q}).$$


Let $E$ be a bundle of rank $n$ and denote by $x_{1},\ldots,
x_{n}$ its Chern roots. The Chern character $ch(E)$ is defined by
the formula:
$$ch(E)=\sum_{i=1}^{n}e^{x_{i}},$$ where $e^{x}=\sum_{n\geq 0}\frac{x^{n}}{n!}$, and $x_{1},\ldots, x_{n}$ are the Chern roots of $E$.
The first few terms are:
$$ch\,(E)=r+c_{1}+\frac{1}{2}(c_{1}^{2}-2c_{2})+\frac{1}{3}(c_{1}^{3}-3c_{1}c_{2}+3c_{3})+\ldots,$$ where $c_{i}=c_{i}(E)$. The $n^{th}$ term
is $\frac{p_{n}}{n!}$, where $p_{n}$ is determined inductively by
Newton's formula (see \cite{Mac}):
$$p_{n}-c_{1}p_{n-1}+c_{2}p_{n-2}-\ldots+(-1)^{n-1}c_{n-1}p_{1}+(-1)^{n}nc_{n}=0.$$

Finally, the Todd class $td\,(E)$ of a bundle with Chern
polynomial $c_{t}(E)=\prod_{i}(1+x_{i}t)$ is defined by
$$td(E)=\prod_{i}\frac{x_{i}}{1-e^{-x_{i}}}.$$

Let us denote by ${\rm{ch}}_{i}$ the $i-$homogeneous part of the
Chern character of a bundle. Then, we have

$$\ch_{0}\,(\mathcal{K}^{\vee}\otimes
\pi_{2}^{*}L)=m=2, $$

$$\ch_{1}\,(\mathcal{K}^{\vee}\otimes
\pi_{2}^{*}L)=t_{1}+\alpha_{1}+\eta\,(d+2a), $$

\begin{align*}\ch_{2}\,(\mathcal{K}^{\vee}\otimes \pi_{2}^{*}L) &=\frac{1}{2}
[t_{1}^{2}+\alpha_{1}^{2}+2\,t_{1}\,\alpha_{1}+2t_{1}\eta\,(d+2a)\\& -2t_{2}-2\alpha_{2}-2u_{1}\eta-2a\eta t_{1}],
\end{align*}

\begin{align*}\ch_{3}\,(\mathcal{K}^{\vee}\otimes \pi_{2}^{*}L) &=
\frac{1}{6}[t^{3}_{1}+\alpha_{1}t_{1}^{2}+3\,(d+2a)\eta\,t_{1}^{2}+2\alpha_{1}
t_{1}^{2}+3\alpha_{1}^{2}t_{1}+\alpha_{1}^{3}
\\& -3t_{1}t_{2}-3t_{1}\alpha_{2}-3\eta u_{1}t_{1}-3\alpha_{1}\alpha_{2}-3\eta (2a+d)t_{2}
+3a\eta\,t_{1}^{2}-3\alpha_{1}t_{2}]
\end{align*}


We observe that $td\,(R_{C,d})=1$ and $td\,(R_{C,d}\times
C)=1+(1-g)\,\eta$, 
putting this together with the computations of the Chern classes
of $\mathcal{K}^{\vee}\otimes \pi_{2}^{*}L$ and the
Grothendieck-Riemann-Roch formula, yields:

\begin{align*}\ch\,(\pi_{1*}(\mathcal{K}^{\vee}\otimes
\pi_{2}^{*}L))=
 \pi_{1*}\left(
(1+(1-g)\eta)\,(\ch\,(\mathcal{K}^{\vee}\otimes
\pi_{2}^{*}L)\right ).\end{align*}

The $i-$homogeneous term of
$\ch(\pi_{1*}(\mathcal{K}^{\vee}\otimes \pi_{2}^{*}L))$ is given
by the formula:
$$
\ch_{i}(\pi_{1*}(\mathcal{K}^{\vee}\otimes
\pi_{2}^{*}L))=(1-g)\,\ch_{i}(\mathcal{K}^{\vee}\otimes
\pi_{2}^{*}L|_{\{q\}\times
R_{C,d}})+{\rm{coeff}}_{\eta}\left(\ch_{i+1}(\mathcal{K}^{\vee}\otimes
\pi_{2}^{*}L)\right).
$$
The first few terms are:
$$\ch_{0} =\rank\,(\pi_{1*}(\mathcal{K}^{\vee}\otimes
\pi_{2}^{*}L))=d+2a+2(1-g),$$
$$\ch_{1}=t_{1}\,(d+2a)+\alpha_{1}\,(d+2a)-at_{1}-a\alpha_{1}-u_{1}+(1-g)\alpha_{1},$$
\begin{align*}
&\ch_{2}=(1-g)[\frac{1}{2}\,\,t_{1}^{2}+\frac{1}{2}\,\alpha_{1}^{2}
+t_{1}\alpha_{1}-t_{2}-\alpha_{2}+3(d+2a)t_{1}^{2}+3(d+2a)\alpha_{1}^{2}] &\\
&+6(d+2a)\alpha_{1}t_{1}-3u_{1}t_{1}-3at_{1}^{2}-3\alpha_{1}^{2}-6a\alpha_{1}t_{1}-3u_{1}\alpha_{1}-3(d+2a)t_{2}\\
&-3(d+2a)\alpha_{2}+3u_{2}+3at_{2}+3\alpha_{2}\\
&\ch_{3}=(1-g)\,[t_{1}^{3}+\alpha_{1}^{3}+3t_{1}\alpha_{1}^{2}+3t_{1}^{2}\alpha_{1}-3t_{1}t_{2}
 -3\alpha_{1}\alpha_{2}-3\alpha_{1}t_{2}+3t_{3}+3\alpha_{3}]
\end{align*}
\begin{center}
$\vdots$
\end{center}

Finally, we get that the Chern classes of
$\pi_{1*}(\mathcal{K}^{\vee}\otimes \pi_{2}^{*}L)$ are given by
the recursive formula:
$$c_{n}(\pi_{1*}(\mathcal{K}^{\vee}\otimes
\pi_{2}^{*}L))=-\sum_{r=1}^{n}\frac{(-1)^{r-1}}{n}r!\ch_{r}(\pi_{1*}(\mathcal{K}^{\vee}\otimes
\pi_{2}^{*}L))\,c_{n-r}(\pi_{1*}(\mathcal{K}^{\vee}\otimes
\pi_{2}^{*}L)).$$


Let $\sigma_{i}$ be the $i-$symmetric function, that is,
$$\sum_{r=0}^{n}\sigma_{r}t^{r}=\prod_{i=1}^{n}(1+x_{i}t),$$ 

\noindent for each $r\geq 1$ the $r'th$ power sum is:
$$p_{r}=\sum x_{i}^{r}=m_{(r)}.$$ The generating function for the
$p_{r}$ is: \begin{equation}\label{eq1} p(t)=\sum_{r\geq
1}p_{r}t^{r-1}=\sum_{i\geq 1}\sum_{r\geq 1}x_{i}^{r}t^{r-1}=
\sum_{i\geq 1}\frac{d}{dt}log\frac{1}{1-x_{i}t}\end{equation}
\begin{equation}\label{eq2}P(t)=\frac{d}{dt}\prod_{i\geq
1}(1-x_{i}t)^{-1}=\frac{d}{dt}log H(t)=\frac{H'(t)}{H(t)}$$
$$P(-t)=\frac{d}{dt}log\, E(t)=\frac{E'(t)}{E(t)}.\end{equation}
From (\ref{eq1}) and (\ref{eq2}), we get that
$$n\sigma_{n}=\sum_{r=1}^{n}(-1)^{r-1}p_{r}\sigma_{n-r}.$$ 

This is a standard formula for symmetric functions (see
\cite{Mac}). The first Chern classes are:

$$c_{1}\,(\pi_{1*}(\mathcal{K}^{\vee}\otimes \pi_{2}^{*}L))=(a+d+1-g)t_{1}+\alpha_{1}-u_{1}+\eta\,(1-g)(d+2a).$$

\begin{align*}
c_{2}\,(\pi_{1*}(\mathcal{K}^{\vee}\otimes
\pi_{2}^{*}L))&=\frac{1}{2}\,c_{1}^{2}-\frac{1}{2}(1-g+d-a)t_{1}^{2}-\frac{1}{2}(1-g)\alpha_{1}^{2}
\\
&-(1-g)(\alpha_{1}+t_{1}+(a+d)\eta t_{1}-t_{2}-\alpha_{2}-u_{1}\eta)\\
& u_{1}t_{1}-\frac{1}{2}(u_{2}+t_{2}).
\end{align*}

\begin{center}
$\vdots$
\end{center}

\begin{coro}
If $C$ is a curve of genus 1, then
$$[R_{C,d,0}]=-c_{1}(\pi_{1*}(\mathcal{K}^{\vee}\otimes
\pi_{2}^{*}L))=-(d+a)t_{1}-\alpha_{1}+u_{1}.$$

\begin{align*}[R_{C,d,-1}]=-c_{2}(\pi_{1*}(\mathcal{K}^{\vee}\otimes
\pi_{2}^{*}L))&=-\frac{1}{2}\,((d+a)t_{1}+\alpha_{1}-u_{1})^{2}
\\&-\frac{1}{2}\,(d-a)t_{1}^{2}+u_{1}t_{1}-\frac{1}{2}(u_{2}+t_{2}).
\end{align*}
\end{coro}

\vspace{0.5cm}
\subsubsection*{Acknowledgments}
I would like to thank Prof. Peter Newstead for comments and
suggestions on both content and presentation of this paper. I also
want to thank Rajagopalan Parthasarathy for pointing out the
reference \cite{Mac}. Finally I would like to thank Prof. J. E.
Andersen and Alex Bene for reading a last version of the paper and
some corrections. This work has been partially supported by The
European Contract Human Potential Programme, Research Training
HPRN-CT-2000-00101 and MPIM grant.

\end{document}